\documentclass[a4paper,12pt]{article}
\usepackage[utf8]{inputenc}
\usepackage{amsmath}
\usepackage{amsfonts}
\usepackage{amssymb}
\usepackage{amsthm}
\usepackage{graphicx}

\usepackage{url}
\usepackage{caption}
\usepackage[numbers,sort&compress]{natbib}

\title{Demonstrating the Applicability of PAINT to Computationally Expensive Real-life Multiobjective Optimization}
\author{Markus Hartikainen and Vesa Ojalehto}
\date{\today}
\begin{document}
\maketitle
\abstract{We demonstrate the applicability of a new PAINT method to speed up iterations of interactive methods in multiobjective optimization. As our test case, we solve a computationally expensive non-linear, five-objective problem of designing and operating a wastewater treatment plant. The PAINT method  interpolates between a given set of Pareto optimal outcomes and constructs a computationally inexpensive mixed integer linear surrogate problem for the original problem. We develop an IND-NIMBUS\textsuperscript{\textregistered} PAINT module to combine the interactive NIMBUS method and the PAINT method and to find a preferred solution to the original problem. With the PAINT method, the solution process with the NIMBUS method take a comparatively short time even though the original problem is computationally expensive.}

\section{Introduction}

In this section, we give background for our study and a brief overview of this paper. First, in Section \ref{sec:paper}, we describe the aim of this study and the structure of this paper. In Section \ref{sec:mo}, we introduce the basic concepts of multiobjective optimization that are used in this paper. In Section \ref{sec:comp_exp}, we consider the main issues related to computationally expensive multiobjective optimization. Finally, in Section \ref{sec:wwtp}, we describe our test case, i.e., the multiobjective problem of designing and operating a  wastewater treatment plant.

\subsection{About this Paper} \label{sec:paper}

In this paper, we demonstrate how the interpolation method PAINT (introduced in \cite{Hartikainen_MCDM2009,Hartikainen_MMOR,Hartikainen_report4}) can be used to speed up the iterations of an interactive method when solving  computationally expensive multiobjective optimization problems. For this, we revisit a computationally expensive five-objective optimization problem from \cite{Sahlstedt2010} that models designing and operating a wastewater treatment plant.

In our study of the wastewater problem, we used the interactive NIMBUS method, as was used in \cite{Sahlstedt2010}. Compared to the previous study of  \cite{Sahlstedt2010}, the iterations of the interactive method NIMBUS were much faster because of a surrogate problem constructed with the PAINT method. In our study, the decision maker was Mr.~Kristian Sahlstedt as in \cite{Sahlstedt2010} and, thus, we were able to ask the decision maker to compare his experiences of using the NIMBUS method to solve the problem with and without the PAINT method. This comparison of the two approaches gave a unique perspective to our study. 

The structure of this paper is as follows: After describing the background of our study in this section, we outline the PAINT approach to solving computationally expensive problems and the PAINT method in Section \ref{sec:PAINT}. We describe the new IND-NIMBUS\textsuperscript{\textregistered} PAINT module (also called the PAINT module in this paper for short) in Section \ref{sec:soft}. In Section \ref{sec:all}, we illustrate how we used the PAINT method to construct a Pareto front approximation  and a mixed integer linear surrogate problem for the wastewater treatment problem. In Section \ref{sec:all}, we also describe our decision maker's involvement in solving the problem with the PAINT module. In Section \ref{sec:analysis}, we further analyze the decision making process of Section \ref{sec:all}. Finally, in Section \ref{sec:conclusions}, we give our conclusions and ideas for further research.

\subsection{Multiobjective Optimization} \label{sec:mo}

Multiobjective optimization concerns simultaneously optimizing multiple conflicting objectives. A general formulation for a multiobjective optimization problem with $k$ objectives is
\begin{equation}
\begin{array}{ll}
\min & (f_1(x),\ldots,f_k(x))\\
\text{s.t. }&x\in S,
\end{array}
\end{equation}
where $f_i$ are the objective functions and $S$ is the feasible set. A vector $x\in S$ is called a (feasible) solution. For these problems, instead of a single optimal solution there typically exist many Pareto optimal solutions. A solution $x\in S$ is said to (Pareto) dominate another solution $y\in S$ if $f_i(x) \leq f_i(y)$ for all $i=1,\ldots,k$ and $f_j(x) < f_j(y)$ for at least one $j\in\{1,\ldots,k\}$. A solution $x^*\in S$ is Pareto optimal, if there does not exist a solution $x\in S$ that dominates it. A vector $z=(f_1(x),\ldots,f_k(x))$ with $x\in S$ is called an outcome, and an outcome is called Pareto optimal if it is given by a Pareto optimal solution. The set of Pareto optimal outcomes is called the Pareto front.

Although many Pareto optimal solutions typically exist, only one has to be chosen for implementation. Distinguishing between Pareto optimal solutions requires preference information about the objectives of the problem. In multiobjective optimization, it is often assumed that there exists a decision maker who is an expert in the application area and who is prepared to answer questions concerning those preferences. In this paper, this whole process of choosing a single solution for implementation is called solving the problem and, when we want to emphasize the decision maker's involvement, it is also referred to as the decision making process. 

The type of information that is asked from the decision maker depends on the multiobjective optimization method that is used to solve the problem. Different types of multiobjective optimization methods (as categorized in \cite{Miettinen1999,Sawaragi1985}) are no-preference methods, a priori methods, a posteriori methods and interactive methods. In no-preference methods the decision maker is not asked any questions. No-preference methods are applicable to problems, where the decision maker is not available or does not want to get involved. In a priori methods, the decision maker is first asked for preference information and then the best solution according to those preferences is found. The difficulty with a priori methods is that the decision maker may find it hard to define preferences without ever seeing feasible or Pareto optimal solutions. In a posteriori methods, a representative set of the Pareto optimal solutions is found from which the decision maker is allowed to choose a preferred one. The difficulty with a posteriori methods is that generating a representative solution set may be time-consuming and choosing a preferred solution from a large set of solutions may be hard (see e.g., \cite{Larichev1992}).

In this paper, we follow the ideology of interactive methods in solving multiobjective optimization problems. In interactive methods, it is understood that any preference information given by the decision maker is only partial and perhaps flawed. Thus, the decision maker is allowed to explore the Pareto optimal solutions by guiding the interactive method. This allows the decision maker to learn about the problem (as argued e.g., in \cite{Miettinen2008}) and find a preferred solution without examining too many solutions. For more information about interactive methods, see e.g., \cite{Miettinen1999,Miettinen2008}. More specifically, in this paper, we use the interactive synchronous NIMBUS method, introduced in \cite{Miettinen1995,Miettinen2000,Miettinen2006}.

\subsection{Solving Computationally Expensive Multiobjective Optimization Problems} \label{sec:comp_exp}

Some multiobjective optimization problems are computationally expensive (see e.g., \cite{Hakanen2011,Hasenjager2005,Aittokoski2008b,Xu2004,Sahlstedt2010}). This may be caused e.g., by the need to use computationally expensive simulations for evaluating the objective functions. Interactive methods have an advantage to a posteriori methods in solving computationally expensive problems, because the decision maker may guide the search in interactive methods and, thus, fewer solutions need to be computed. There is, however, a drawback. When using interactive methods, the decision maker has to wait while new solutions are computed with respect to his/her updated preferences. For computationally expensive problems, this may take a long time, which may be frustrating for the decision maker (as argued e.g., in \cite{Korhonen1996}).

In order to compute new solutions faster within the interactive method, one can use approximation. Two different approximation schemes can be identified: approximating the objective functions and approximating the Pareto front. The objective functions may be approximated e.g., with meta-models like the response surface methodologies, Support Vector Machines or Radial Basis Functions (see e.g., \cite{Nakayama2009}). These have been used in multiobjective optimization e.g., in \cite{Nakayama2009,Wilson2001}. This is not, however, a straightforward task, because as the number of decision variables and objectives increases, the approximation itself becomes a very computationally expensive task. Another approach is approximating the Pareto front. Pareto front approximations can be found e.g., in \cite{Ackermann2007,Bezerkin2006,EskelinenMiettinenKlamrothHakanen2008,Lotov2004,Monz2008,RuzikaWiecek2005,Yapicioglu2011}, where \cite{Ackermann2007,EskelinenMiettinenKlamrothHakanen2008,Lotov2004,Monz2008} include decision making aspects connected to these. Note that in this paper, we distinguish between a Pareto front representation (a discrete set of Pareto optimal outcomes) and a Pareto front approximation (something more approximate that possibly contains vectors that are not outcomes of the problem, but merely approximate them). 

In this paper, we use the Pareto front approximation approach introduced in \cite{Hartikainen_MCDM2009,Hartikainen_MMOR,Hartikainen_report4}. In those papers, a new Pareto front approximation method PAINT is introduced and details on decision making with the produced approximation are covered. The PAINT method uses a novel way to integrate the knowledge about Pareto dominance into the approximation. The PAINT method interpolates between a given set of Pareto optimal outcomes to construct a Pareto front approximation. The approach differs from the other approaches for decision making with Pareto front approximations (mentioned above) because it is able to approximate also nonconvex Pareto fronts. Furthermore, the Pareto front approximation constructed with PAINT implies a multiobjective mixed integer linear surrogate problem (for the original problem) that can be solved with any interactive method. The other approaches are either applicable only to convex multiobjective optimization problems or use only a custom-made procedure for choosing a preferred point on the approximation. Further details on the PAINT method are covered in Section \ref{sec:PAINT}.

\subsection{Designing and Operating a Wastewater\\*Treatment Plant} \label{sec:wwtp}

Designing and operating a wastewater treatment plant is a complex problem with many conflicting criteria that have to be considered at the same time. In this paper, we consider a plant using so-called activated sludge process, which is globally the most common method of wastewater treatment. We model the problem as a five-objective optimization problem, which was previously studied also in \cite{Sahlstedt2010}. The five-objective problem is an extension of the three-objective problem treated in \cite{Hakanen2011}. The approach of this paper differs from the approach of \cite{Sahlstedt2010} because we use the PAINT method to approximate the Pareto front and to construct a surrogate problem for the original problem. In this way, the time that the decision maker has to wait while using an interactive method becomes shorter.

Figure \ref{fig:schematic} shows the schematic layout of the wastewater treatment plant that was designed in \cite{Sahlstedt2010}. The wastewater treatment begins with grit removal. After the grit removal, solids are separated by a gravitational settling. Raw and mixed sludge removed from the primary settlers is fermented in a separate reactor and partly recycled back to the water line to provide readily biodegradable carbon source for denitrification. The bioreactor consists of four anoxic zones, three aerobic zones and one deoxygenation zone. Nitrate-rich activated sludge is recycled from zone 8 of the bioreactor to zone 1. Return sludge and primary effluent are directed to zone 1. Methanol is injected to zone 2 to support denitrification. Excess sludge is pumped from zone 8 of the bioreactor to the beginning of the water process, from which it is removed in the primary settlers together with raw sludge. Raw and mixed sludge is thickened gravitationally into approximately 4.5\% total solids prior to anaerobic digestion. Anaerobic digestion produces biogas and the produced biogas can be converted into electrical or thermal energy. The digested sludge is dewatered by centrifuges into approximately 28\% total solids. The reject water from sludge treatment is pumped to the beginning of the plant. The wastewater treatment process is simulated with the commercial GPS-X simulator (see \cite{GPS-X_homepage}) and the model is based on the findings of Pöyry Engineering Ltd. For more information about the wastewater treatment plants using activated sludge process, see e.g., \cite{Hakanen2011,Philips2009,Sahlstedt2010}.

The objectives of the optimization problem are the amount of nitrogen in the effluent ($g/m^3$, grams per a cubic meter of effluent), aeration power consumption in the activated sludge process ($kW$), chemical consumption ($g/m^3$, grams per a cubic meter of effluent), excess sludge production ($kg/d$, kilograms per day) and biogas production ($m^3/d$, cubic meters per day). The first one is the main goal of activated sludge process and the four others are connected to the operational costs. This multiobjective optimization problem allows the simultaneous consideration of the performance of the plant (through the nitrogen removal rate) and different aspects of the operational costs. Naturally, the last objective is maximized and the others are minimized. Decision variables of the problem are the percentage of inflow pumped to fermentation, the amount of excess sludge removed, the dissolved oxygen setpoint in the last aerobic zone and the methanol dose. Thus, the methanol dose is both a decision variable and an objective.

\begin{center}
    \includegraphics[width=0.80\textwidth]{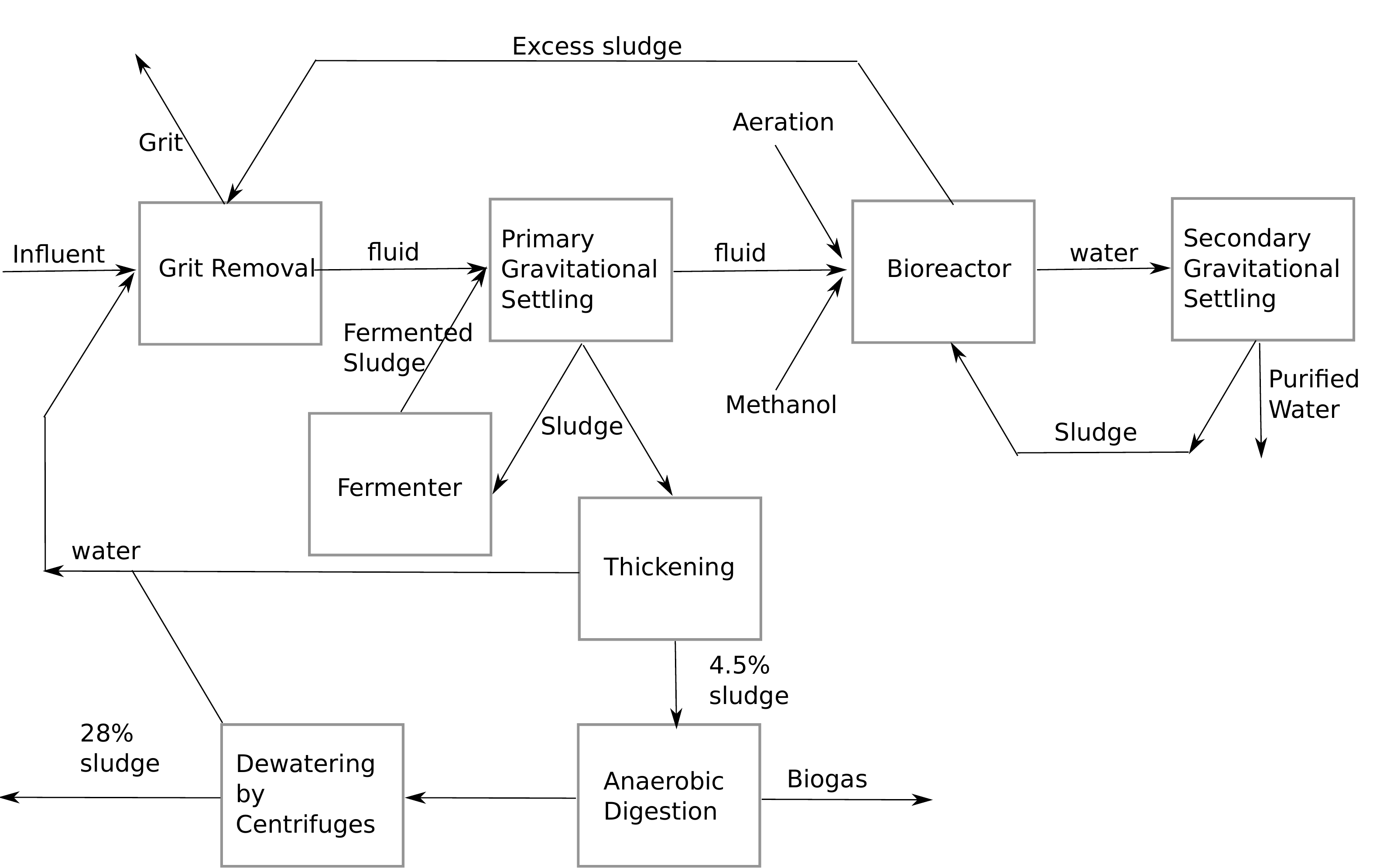}
    \captionof{figure}{A schematic layout of the wastewater treatment plant}
    \label{fig:schematic}
\end{center}

Each simulation of the wastewater treatment plant of \cite{Sahlstedt2010} took about $11$ seconds on the GPS-X simulator. This made the problem computationally expensive. In addition, one could notice from the Pareto optimal outcomes computed for the problem that the problem is nonconvex. During the analysis in \cite{Sahlstedt2010}, $200$ simulations were run to optimize the scalarizations (i.e., single objective optimization problems, whose optimal solutions are Pareto optimal solutions to the multiobjective optimization problem) given by the interactive NIMBUS method that was used to solve the problem. This means that each iteration of the interactive method took more than half an hour. Even though interesting solutions to the problem were found in \cite{Sahlstedt2010}, the computational time of iterations was an inconvenience to the decision maker (according to personal communications with the authors of \cite{Sahlstedt2010}). This means that there was room for improvement using the PAINT method.

\section{The PAINT Approach to Solving Computationally Expensive Problems} \label{sec:PAINT}

In this section, we describe the PAINT approach to solving computationally expensive problems. The applicability of the PAINT approach is then demonstrated in Section \ref{sec:all} by solving a computationally expensive multiobjective optimization problem of wastewater treatment plant design and operation.

The PAINT approach is based on the Pareto front approximation constructed by the PAINT method. The PAINT method was proposed in \cite{Hartikainen_report4}, and it is based on the concept of an inherently nondominated Pareto front approximation introduced in \cite{Hartikainen_MCDM2009} and the mathematical concepts of \cite{Hartikainen_MMOR}. The PAINT method interpolates between a given set of Pareto optimal outcomes in a way that the interpolants neither dominate nor are dominated by the set of given Pareto optimal outcomes and, in addition, they are not dominated by each other (i.e., the interpolation is an inherently nondominated Pareto front approximation, as defined in \cite{Hartikainen_MCDM2009}). In this paper, a vector on the Pareto front approximation is called an approximate (Pareto optimal) outcome.

The general functionality of the PAINT method is as follows: The PAINT method first constructs the Delaunay triangulation of the given set of Pareto optimal outcomes and then chooses the appropriate polytopes from it to the Pareto front approximation. In this paper, this is realized with the Octave-based (see \cite{Octave_homepage,Eaton2002}) implementation that was developed during the research of \cite{Hartikainen_report4}. 

The Pareto front approximation constructed with the PAINT method implies a computationally inexpensive mixed integer linear surrogate problem for the original problem, as described in \cite{Hartikainen_report4}. The Pareto front of the surrogate problem is exactly the Pareto front approximation and, thus, a preferred solution to the surrogate problem implies a preferred vector on the Pareto front approximation, which is also called a preferred approximate outcome in this paper. The algorithm of the PAINT method and more exact details can be found in \cite{Hartikainen_report4}. 

\begin{figure}
\begin{center}
    \includegraphics[width=0.60\textwidth]{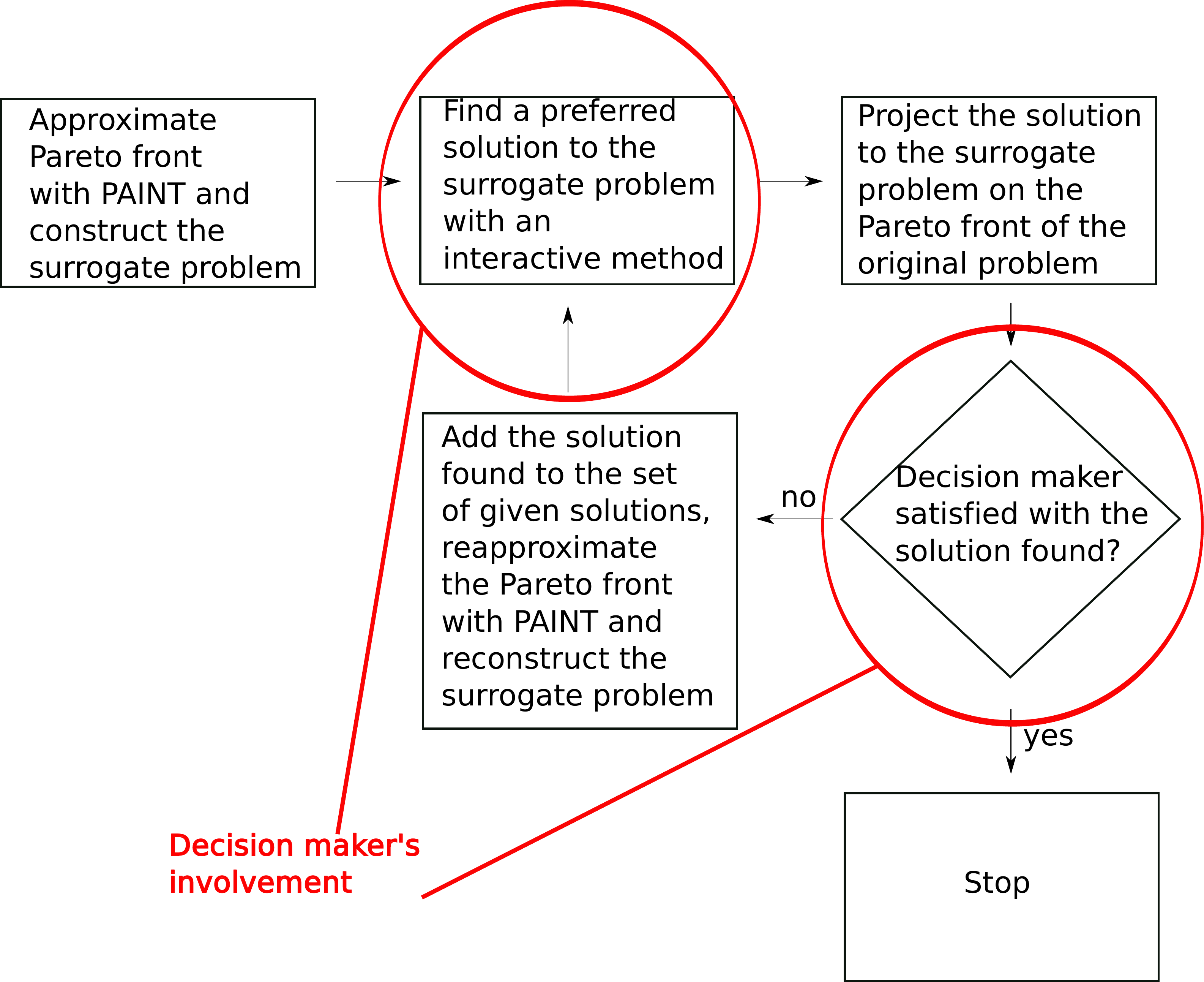}
\end{center}
\caption{A flowchart of the decision making process}
\label{fig:DM}
\end{figure}

Decision making in the PAINT approach is described in Figure \ref{fig:DM}. In the PAINT approach to solving computationally expensive problems, we assume that there exists a set of Pareto optimal solutions to the computationally expensive problem. This set may have been generated with any a posteriori method. The set of the related outcomes is inputted into the PAINT method. PAINT then interpolates between the set of given Pareto optimal outcomes and outputs the interpolation that implies a mixed integer linear surrogate problem for the original problem.  

After the mixed integer linear surrogate problem has been formulated, the decision maker gets involved and uses an interactive method of his/her choice to find a preferred solution to the surrogate problem. The outcomes given by Pareto optimal solutions to the surrogate problem are vectors on the Pareto front approximation, which is in the same space as the original Pareto front. Thus, the decision maker is able to give his/her preferences on them. The preferred approximate outcome is projected on the actual Pareto front of the original problem by solving achievement scalarizing problem (see \cite{Wierzbicki1986,Wierzbicki1980b}) with the approximate outcome as a reference point. More details on the projection can be found in \cite{Hartikainen_report4}. Projecting the solution may take time, depending on the computational costs of the problem. If the problem is very computationally expensive, the projection can be done without the involvement of the decision maker. 

The projection of the preferred approximate outcome (i.e., a Pareto optimal solution to the original problem) is shown to the decision maker and, if he/she is satisfied, the decision making process stops, because a preferred solution has been found. If the decision maker is not satisfied, it is possible to update the Pareto front approximation by adding the new Pareto optimal outcome to the given set of Pareto optimal outcomes and by recomputing the approximation with the PAINT method. This yields a more accurate approximation and we can again use an interactive method to find a preferred solution to the new (more accurate) surrogate problem. This process can be repeated as many times as necessary.       

The PAINT method is a powerful tool as it can interpolate between any given set of Pareto optimal outcomes, i.e., the way that the outcomes have been generated does not affect the functionality of the method. In addition, since it is based on the concept of inherent nondominance (see \cite{Hartikainen_MCDM2009}), it will not provide interpolants that would mislead the decision maker. Finally, the mixed integer linear surrogate problem implied by the approximation allows one to use any interactive method for finding a preferred approximate outcome on the Pareto front approximation. 

The PAINT method has a couple of shortcomings, already noted in \cite{Hartikainen_report4}. First, the PAINT method does not provide any information about the preimage of the Pareto front approximation in the decision space. This means that the decision maker has to project the approximate outcome (i.e., the solution to the surrogate problem) on the Pareto front of the original problem in order to find out the values of the decision variables. Second, the PAINT method cannot detect any disconnectedness in the Pareto front, but always interpolates between the outcomes whenever the interpolation is inherently nondominated. Thus, the approximation might be inaccurate if e.g., the decision space is disconnected or the objective functions are highly nonconvex.

\section{A New IND-NIM\-BUS\textsuperscript{\textregistered} PAINT Module}  \label{sec:soft}

IND-NIMBUS\textsuperscript{\textregistered} (see \cite{Miettinen2006b}) is a multi-platform desktop software framework, currently available for Windows and Linux operating systems, intended to provide a flexible tool-set for implementation of multiobjective optimization methods. So far, the IND-NIMBUS framework has been used to implement the synchronous NIMBUS \cite{Miettinen1995,Miettinen2000,Miettinen2006} and the Pareto Navigator \cite{EskelinenMiettinenKlamrothHakanen2008} methods. The IND-NIMBUS\textsuperscript{\textregistered} software can be connected to an external sources that model the problem, such as the GPS-X simulator used for modeling the wastewater treatment plant. For this paper, the IND-NIMBUS\textsuperscript{\textregistered} software framework has been used to develop a so-called IND-NIMBUS\textsuperscript{\textregistered} PAINT module that combines the PAINT and NIMBUS methods for computationally expensive multiobjective optimization. The PAINT module implements most of the functionalities described in Figure \ref{fig:DM}.

The synchronous NIMBUS method \cite{Miettinen1995,Miettinen2000,Miettinen2006} is an interactive multiobjective optimization method.  The NIMBUS method uses classification of objectives as the preference information. Given a Pareto optimal solution to the multiobjective optimization problem, the decision maker can classify the objectives into classes $I^<$, $I^{\leq}$, $I^=$, $I^{\geq}$ and $I^{<>}$ defined, respectively, as classes of objective functions that the decision maker wants to improve as much as possible, wants to improve to a given aspiration level $z_i$, allows to remain unchanged, allows to deteriorate until a given bound $\epsilon_i$ and allows to change freely for a while. This preference information is converted into several different single objective subproblems with the help of different scalarization functions as proposed in \cite{Miettinen2006}. These subproblems are solved to generate different Pareto optimal solutions, which are shown to the decision maker who can then see how well the desired preferences could be attained. The decision maker can choose any of these solutions as the starting point of the next iteration, i.e., classification. This iterative procedure can either start with a solution given by the decision maker or from a so-called neutral compromise solution and it is repeated until the decision maker is satisfied with the solution at hand. Further information about the synchronous NIMBUS with other means to direct the search process is given in \cite{Miettinen2006}.

The NIMBUS method has been successfully applied to shape design of ultrasonic transducers \cite{HeiMieNie06}, designing a paper machine headbox \cite{HamMieTarToi03}, optimal control in continuous casting of steel \cite{Mie07}, separation of glucose and fructose \cite{HakKawMieBie07}, intensity modulated radiotherapy treatment planning \cite{RuoBomMieTer09}, brachytherapy \cite{Ruotsalainen2010} and optimizing heat exchanger network synthesis \cite{Laukkanen2010}, among others. In addition, it uses  classification of objectives that has been found cognitively just \cite{Larichev1992}. These facts make the NIMBUS method an ideal choice as the interactive method for solving the PAINT surrogate problem of the wastewater treatment plant model.

\begin{figure}[h!]
\begin{center}
    \includegraphics[width = 0.9\textwidth]{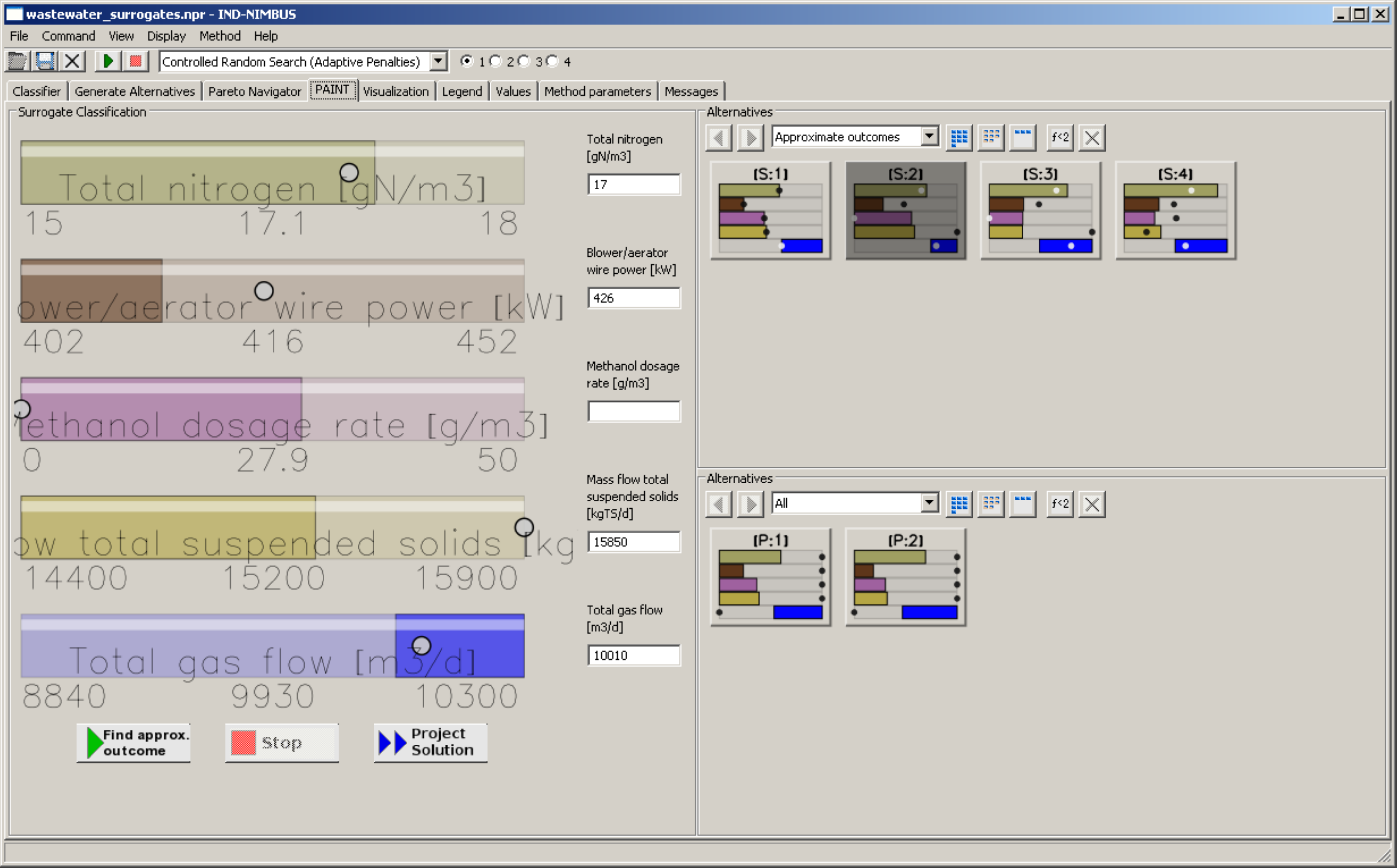}
\end{center}
   \caption{A screen shot of the IND-NIMBUS\textsuperscript{\textregistered} PAINT module}
    \label{fig:soft}
\end{figure}

Figure \ref{fig:soft} shows the screen shot of the IND-NIMBUS\textsuperscript{\textregistered} PAINT module. In the PAINT module, the decision maker can give his/her preferences concerning the surrogate problem by classification of the objective functions. The given classification information is used to formulate a single objective subproblem of the surrogate problem. The subproblem is modeled using  the Optimizing Programming Language (OPL, see \cite{Henterryck1999}), and this (mixed integer linear model) is solved using CPLEX (see \cite{Cplex_homepage}). An optimal solution to the subproblem gives a new approximate Pareto optimal outcome, corresponding to the preferences given by the decision maker. This approximate outcome is shown to the decision maker. If the decision maker so wishes, he reclassify the objectives of the new approximate outcome which yields another approximate outcome.

As described in Section \ref{sec:PAINT}, approximate Pareto optimal outcomes can be projected on the Pareto front of the original problem using the PAINT module (using the Project Solution button near the bottom of the screen).  The projection of the approximate outcome, that is, a Pareto optimal solution to the original problem is shown to the decision make. As mentioned, for a computationally expensive problem this may take time, but fortunately projecting an approximate outcome can be done without the involvement of the decision maker. 

The approximate Pareto optimal outcomes and the actual Pareto optimal solutions that have been found during the decision making process are visualized on the right side of the PAINT module. The decision maker can choose any of the approximate Pareto optimal outcomes as the starting point of the next NIMBUS iteration (i.e., as the basis for a new classification of objectives). The process stops when the decision maker has found a preferred solution to the original problem through projecting an approximate Pareto optimal outcome.

In the current version of the PAINT module, only one of the scalarizations of the synchronous NIMBUS method (i.e., the achievement scalarizing subproblem) has been implemented. That is, unlike in the synchronous NIMBUS method, the decision maker can see only one approximate Pareto optimal outcome for given preferences. It should also be noted that any solver capable of solving the surrogate subproblem (e.g., GLPK, see \cite{glpk_homepage}) could be used instead of CPLEX. 

The current version of the PAINT module does not implement the construction or updating of the surrogate problem.  If one wishes to update the surrogate problem using Pareto optimal outcomes obtained by e.g., projection, the decision making process must be stopped, and the surrogate problem must be manually updated using Octave. In future versions, updating the Pareto front approximation should be implemented under a third button in the PAINT module that would then automatically update the approximation and the surrogate problem.

\section{Solving the Wastewater Treatment Case} \label{sec:all}

In this section, we demonstrate how the PAINT method and the IND-NIMBUS\textsuperscript{\textregistered} PAINT module were used to solve the wastewater treatment problem, described in Section \ref{sec:wwtp}. First, in Section \ref{sec:analysis_without_DM}, we describe the construction of the Pareto front approximation with PAINT before the involvement of the decision maker. Then, in Section \ref{sec:analysis_with_DM}, we describe how the decision maker used the PAINT module to solve the wastewater treatment problem.

\subsection{Pre-Decision Making phase} \label{sec:analysis_without_DM}

First, a set of $200$ mutually nondominated solutions to the wastewater problem was found with the evolutionary UPS-EMO algorithm (introduced in \cite{Aittokoski2010}) and the GPS-X simulator. To study the optimality of these solutions, each one was locally improved using an achievement scalarizing problem \cite{Wierzbicki1980b,Wierzbicki1986}, which was optimized with Matlab fmincon-function with finite differences approximated gradients. This resulted in $195$ mutually nondominated solutions. The maximal improvement in the values of the achievement scalarizing problem was at most $3\%$ so the local improvement did not cause much change. This built our confidence that the final solutions were close to Pareto optimal. We took the outcomes given by these solutions as the set of given Pareto optimal outcomes for the PAINT method. The whole process of producing this set took about three days on a standard laptop. 

After this, we computed a Pareto front approximation based on the given set of Pareto optimal outcomes with the PAINT method (see Section \ref{sec:PAINT}). The PAINT method chose $4272$ polytopes for interpolation in the Pareto front. In order to reduce the computational complexity of the implied mixed integer linear surrogate problem, we removed polytopes that were subsets of larger polytopes from the approximation. This resulted in a collection of $608$ polytopes whose union covered the same space in $\mathbb R^5$ as that of the larger collection. In addition, all sets of vertices of the polytopes in the collection were affinely independent and, thus, the number of vertices of all the polytopes was five or less, as shown in \cite{Hartikainen_MMOR}. Using the PAINT method to construct the Pareto front approximation took approximately 19 hours on Intel\textsuperscript{\textregistered} Xeon\textsuperscript{\textregistered} E5410 CPU. 

The mixed integer linear surrogate problem implied by the smaller collection was equivalent to that implied by the larger collection, but it was computationally less expensive. As described in \cite{Hartikainen_report4}, the surrogate problem could be written as
\begin{equation} \label{eq:surrogate}
\begin{array}{ll} \min &(z_1,\ldots,z_5)\\
    \text{s.t. } &\sum_{j=1}^{608} \sum_{l=1}^5 \lambda_{j,l} = 1\\
                 &\sum_{l=1}^5 \lambda_{j,l} \leq y_j, \text{ for all }j=1,\ldots,{608}\\
                 &\sum_{j=1}^{608} y_j = 1\\
    \text{where }&\lambda\in[0,1]^{608\times 5}\\
            &y\in \{0,1\}^{608}\\
             &z_i = \sum_{j=1}^{608} \sum_{l=1}^5 \lambda_{j,l}p^{A_{l,j}}_i\text{ for all }i=1,\ldots,5,
    \end{array}
\end{equation}
where each row of the matrix $A\in \mathbb R^{608\times 5}$ contained the indices of the vertices of one polytope in the smaller collection of polytopes. The component $\lambda_{j,l}$ of the matrix variable $\lambda\in\mathbb R^{608\times 5}$ was for all $j=1,\ldots,608$ and $l=1,\ldots,5$ the coefficient of the vertex $l$ of the polytope given by row $j$ in the matrix $A$. The variable $y$ determined which of the rows of the matrix $\lambda$ was nonzero. By the third constraint, only one row in the matrix $\lambda$ had nonzero elements.

Problem \eqref{eq:surrogate} had $608\times5=3040$ continuous variables and $608$ binary variables. CPLEX was able to solve e.g., an achievement scalarizing problem for the surrogate problem in less than a second. This was a tremendous improvement to solving a scalarization of the original problem, which took about half an hour with the Controlled Random Search algorithm (see \cite{Price1983}), which is implemented in the IND-NIMBUS software.

\subsection{Decision Making Phase} \label{sec:analysis_with_DM}

Using the PAINT module, our decision maker (Mr. Kristian Sahlstedt from Pöyry Environment Ltd) was able to examine the approximate outcomes and to project any of them on the Pareto front of the original problem. This entire decision making process was done within a couple of hours and the decision maker's involvement was only about an hour, which could not have been possible by merely using the original computationally expensive problem.  

Before the decision maker started using the PAINT module, we gave him a brief overview of the methods from the user's perspective. We told him that a set of Pareto optimal outcomes has been computed and that a new PAINT method has been used to interpolate between those outcomes. We also informed him that the outcomes given by PAINT are only approximate Pareto optimal outcomes and, thus, more computation has to be done to find the closest real Pareto optimal outcome. In addition, we told him that the PAINT method does not unfortunately provide any information about the decision variables and those values can only be known after the real Pareto optimal solution is found. Since our decision maker had previous experiences with the NIMBUS method, all of this was very clear to him. In addition, he did not find any of this too inconvenient.

\begin{table}
\begin{center}
\begin{tabular}{||p{0.05\textwidth}|p{0.13\textwidth}|p{0.13\textwidth}|p{0.13\textwidth}|p{0.13\textwidth}|p{0.13\textwidth}||}
\hline
Out\-come &Amount of Nitrogen [$gN/m^3$] (min) 	&Aeration power [$kW$] (min) 	&Chemical consumption [$g/m^3$] (min)  &Excess sludge [$kg/d$] (min) 	&Biogas production [$m^3/d$] (max)\\
\hline
$s^1$	&16.67				&412.2			&21.89				&15060				&9731\\
$s^2$	&17.13				&416.3			&27.86				&15250				&9935\\
$s^3$	&17.30				&419.0			&16.27				&14870				&9560\\
$s^4$	&17.74				&414.6			&14.41				&14910				&9571\\
$p^1$	&16.80				&414.1	 		&18.24				&14960				&9626\\
$p^2$	&17.10				&411.6			&15.10				&14860				&9529\\
\hline
\end{tabular}
\caption{The approximate and actual Pareto optimal outcomes inspected by the decision maker with the PAINT module}
\label{tab:sol}
\end{center}
\end{table}

Table \ref{tab:sol} shows the approximate Pareto optimal outcomes generated (approximate outcomes $s^1,\ldots,s^4$) and the outcomes given by actual Pareto optimal solutions to the wastewater treatment problem (outcomes $p^1,p^2$) that were inspected by the decision maker. The decision making process started from the approximate outcome $s^1$ in Table \ref{tab:sol}. The outcome $s^1$ was given by the neutral compromise solution to the surrogate problem.

The decision maker wanted to see further (approximate) Pareto optimal outcomes. After a classification of objectives in the PAINT module, the optimal solution to the new subproblem for the surrogate problem gave the approximate Pareto optimal outcome $s^2$. The approximate outcome $s^2$ has more biogas production than the approximate outcome $s^1$, but is worse in all the other objectives. Especially, the chemical consumption is very large. Let us emphasize that finding the approximate outcome $s^2$ was especially smooth, since the mixed integer linear problem was computationally inexpensive.

Because the decision maker was not completely satisfied with the approximate outcome $s^2$, he decided to continue and find another approximate outcome. This yielded the approximate outcome $s^3$, which has much lower chemical consumption and a slightly lower excess sludge production than both approximate outcomes $s^1$ and $s^2$. Unfortunately, the approximate outcome $s^3$ is worse than both approximate outcomes $s^1$ and $s^2$ in all the other objectives.

A classification of objectives of the approximate outcome $s^3$ and solving the new subproblem for the surrogate problem led to the approximate outcome $s^4$. This new approximate outcome is roughly the same as the approximate outcome $s^3$ in both excess sludge production and biogas production. However, the amount of nitrogen for the approximate outcome $s^4$ is slightly higher than for the approximate outcome $s^3$, but this is compensated by chemical consumption and aeration power that are considerably smaller.

After having inspected the four approximate Pareto optimal outcomes, the decision maker felt that he had learned enough about the surrogate problem. First, the decision maker decided to project the approximate outcome $s^1$ on the Pareto front of the original wastewater treatment problem. The projection of the approximate outcome $s^1$ (i.e., solution to the achievement scalarizing problem with the approximate outcome $s^1$ as the reference  point) took a little over half an hour. The projection was done using the GPS-X simulator and the Controlled Random Search algorithm. The projection was outcome $p^1$ in Table \ref{tab:sol}. According to our decision maker's assessment, the outcome $p^1$ was fairly close to the approximate outcome $s^1$ in all objectives. However, he felt that there might still be more preferred solutions to the problem.

Because the approximate outcome $s^1$ and the actual Pareto optimal outcome $p^1$ were close to each other, it was decided not to update the Pareto front approximation. Instead, the decision maker wanted to project the approximate outcome $s^4$ and obtained the outcome $p^2$ in Table \ref{tab:sol}. The Pareto optimal outcome $p^2$ has slightly higher amount of nitrogen in the effluent than the Pareto optimal outcome $p^1$, but it has considerably lower aeration power and the amount of chemical consumption. The Pareto optimal outcome $p^2$ was very preferred by the decision maker and he chose it as the final solution to the problem. 

The final solution (including a way to design and operate a wastewater treatment plant) will be further inspected by more accurate simulators before implementing. However, it will act as a guideline for the design of the wastewater treatment plant. 

\section{Analyzing the Decision Making Process with PAINT} \label{sec:analysis}

We filmed our decision maker Mr.~Sahlstedt during the decision making process and asked him some additional questions regarding the usability of the methods. The purpose of the video was to reveal any issues that he might have had while using the PAINT module and to find out whether any aspects of the PAINT method that were hard to understand. 

When analyzing the video, it seems that the key point in the usability of the PAINT method and the PAINT module is informing the decision maker about the approximate nature of the method. In order to make the PAINT method more usable, it would be a good idea to produce an introductory video introducing the key points of the method. Since our decision maker was already familiar with the IND-NIMBUS\textsuperscript{\textregistered} software, there was no need to introduce it. This may not always be the case and, thus, the video should also include a short introduction to the NIMBUS method and the IND-NIMBUS\textsuperscript{\textregistered} software.

The decision maker thought that the PAINT method and the PAINT module were easy and intuitive to use. In addition, he said that the PAINT method provided a definite improvement to merely using the NIMBUS method because of the faster computational times (a couple of seconds against half an hour) between the iterations of the interactive method. He thought that neither the approximate nature nor the fact that the preferences had to be based only on the objective function values were big drawbacks. In addition, no approximate outcomes that the decision maker would have assessed implausible were found during the decision making process.

For both projected approximate outcomes $s^1$ and $s^4$, our decision maker assessed that the actual Pareto optimal solutions were close enough, taking into account the uncertainties in the model, and there was no need to recompute the Pareto front approximation. Thus, the inability of the PAINT module in reconstructing the Pareto front approximation (as mentioned in Section \ref{sec:soft}) was not an issue.

Except for correcting a couple of minor bugs from the software, the decision maker did not offer any improvements. He did, however, agree with us that the PAINT method should be able to detect disconnectedness in the Pareto front and that the decision variables should be somehow approximated, too. However, detecting disconnectedness was not an issue on this occasion, because all the approximate Pareto optimal outcomes found seemed plausible to the decision maker and the approximate Pareto optimal outcomes that were projected were rather close to their projections. Finally, since one of the decision variables was also an objective (i.e., the methanol dose), that was approximated in our problem, although the PAINT method does not in general do this.

\section{Conclusions} \label{sec:conclusions}

In this paper, we described how the PAINT method can be used to speed up iterations of interactive methods when solving computationally expensive multiobjective optimization problems. The PAINT method was used to construct a Pareto front approximation that then implied a mixed integer linear surrogate problem for the original problem. As our case problem, we studied a five-objective optimization problem of designing and operating a wastewater treatment plant. In addition, we introduced a new IND-NIMBUS\textsuperscript{\textregistered} PAINT module that combines the PAINT method and the interactive NIMBUS method. 

The PAINT method and the PAINT module worked well in this problem. The decision maker found it easy and intuitive to use the interactive NIMBUS method to find a preferred approximate outcome on the Pareto front approximation. The low computational cost of using the interactive method with the surrogate problem was a definite improvement to using interactive method directly to solve the computationally expensive wastewater treatment problem.

The experimental design in this paper was unique: Because our decision maker had already used the IND-NIMBUS\textsuperscript{\textregistered} to study the same wastewater treatment problem, he was able to compare the experiences of using the IND-NIMBUS software with and without the PAINT method. According to the decision maker's opinion, the PAINT method provided a significant improvement. This implies that the PAINT method should be also applicable to other computationally expensive problems.

The IND-NIMBUS\textsuperscript{\textregistered} PAINT module is still in the development phase and, thus, it lacks some essential functionality (like the implementations of the other scalarizations of the synchronous NIMBUS method) and, also, it still has some bugs. If there had been no bugs in the software, the investigation of the problem with the PAINT module would have been even more fluent. Further effort has to be put in correcting these bugs.

The PAINT method requires one to use additional methods and software to generate the given set of Pareto optimal outcomes and to solve the mixed integer linear surrogate problem. In this paper, we used the UPS-EMO algorithm to generate the Pareto optimal outcomes and the IND-NIMBUS software with the PAINT module and the CPLEX solver to solve the surrogate problem. In future research, other applicable methods and software can be also used together with the PAINT method. 

\section{Acknowledgements}

The authors of this paper were financially supported by the Academy of Finland (grant number $128495$). The authors are indebted to Mr.~Kristian Sahlstedt from Pöyry Environment Ltd for his expertise on the problem and for acting as the decision maker. In addition, the authors wish to thank Drs.~Jussi Hakanen and Timo Aittokoski from the Research Group in Industrial Optimization at the University of Jyväskylä, Department of Mathematical Information Technology for helping to solve the problem. During the writing of this paper, Prof.~Miettinen provided her helpful comments that helped the authors improve this paper.

\bibliography{mybib,references}

\begin{thebibliography}{10}

\bibitem{Ackermann2007}
H.~Ackermann, A.~Newman, H.~R\"oglin, and B.~V\"ocking.
\newblock {Decision-making Based on Approximate and Smoothed Pareto Curves}.
\newblock {\em Theoretical Computer Science}, 378:253--270, 2007.

\bibitem{Aittokoski2008b}
T.~Aittokoski and K.~Miettinen.
\newblock {Cost Effective Simulation-Based Multiobjective Optimization in the
  Performance of an Internal Combustion Engine}.
\newblock {\em Engineering Optimization}, 40:593--612, 2008.

\bibitem{Aittokoski2010}
T.~Aittokoski and K.~Miettinen.
\newblock {Efficient Evolutionary Approach to Approximate the Pareto-Optimal
  Set in Multiobjective Optimization, UPS-EMOA}.
\newblock {\em Optimization Methods and Software}, 25:841 -- 858, 2010.

\bibitem{Bezerkin2006}
V.~E. Bezerkin, G.~K. Kamenev, and A.~V. Lotov.
\newblock {Hybrid Adaptive Methods for Approximating a Nonconvex
  Multidimensional Pareto Frontier}.
\newblock {\em Computational Mathematics and Mathematical Physics},
  46:1918--1931, 2006.

\bibitem{Eaton2002}
J.~W. Eaton.
\newblock {\em {GNU Octave Manual}}.
\newblock Network Theory Limited, 2002.

\bibitem{EskelinenMiettinenKlamrothHakanen2008}
P.~Eskelinen, K.~Miettinen, K.~Klamroth, and J.~Hakanen.
\newblock {Pareto Navigator for Interactive Nonlinear Multiobjective
  Optimization}.
\newblock {\em OR Spectrum}, 32:211--227, 2010.

\bibitem{glpk_homepage}
{GLPK (GNU Linear Programming Kit) Home Page}, Last accessed: \today.
\newblock url = \url{http://www.gnu.org/s/glpk/}.

\bibitem{GPS-X_homepage}
{GPS-X Simulator Home Page}, Last accessed: \today.
\newblock url =
  \url{http://www.technotrade.com.pk/20/GPS-X_Waste_Water_Treatment_Simulation%
_Software/}.

\bibitem{HakKawMieBie07}
J.~Hakanen, Y.~Kawajiri, K.~Miettinen, and L.T. Biegler.
\newblock {Interactive Multi-Objective Optimization for Simulated Moving Bed
  Processes}.
\newblock {\em Control and Cybernetics}, 36(2):282--320, 2007.

\bibitem{Hakanen2011}
J.~Hakanen, K.~Miettinen, and K.~Sahlstedt.
\newblock {Wastewater Treatment: New Insight Provided by Interactive
  Multiobjective Optimization}.
\newblock {\em Decision Support Systems}, 51:328--337, 2011.

\bibitem{HamMieTarToi03}
J.H. H\"am\"al\"ainen, K.~Miettinen, P.~Tarvainen, and J.~Toivanen.
\newblock {Interactive Solution Approach to a Multiobjective Optimization
  Problem in Paper Machine Headbox Design}.
\newblock {\em Journal of Optimization Theory and Applications}, 116:265--281,
  2003.

\bibitem{Hartikainen_MMOR}
M.~Hartikainen, K.~Miettinen, and M.~M. Wiecek.
\newblock {Constructing a Pareto Front Approximation for Decision Making}.
\newblock {\em Mathematical Methods of Operations Research}, 73:209--234, 2011.

\bibitem{Hartikainen_report4}
M.~Hartikainen, K.~Miettinen, and M.~M. Wiecek.
\newblock {PAINT: An Interpolation Method for Computationally Expensive
  Multiobjective Optimization Problems}, 2011.
\newblock University of Jyv\"askyl\"a, Reports of the Department of
  Mathematical Information Technology, Series B. Scientific Computing, Number
  B3/2011.

\bibitem{Hartikainen_MCDM2009}
M.~Hartikainen, K.~Miettinen, and M.~M. Wiecek.
\newblock {Pareto Front Approximations for Decision Making with Inherent
  Non-dominance}.
\newblock In Y.~Shi, S.~Wang, G.~Kou, and J.~Wallenius, editors, {\em {New
  State of MCDM in the 21st Century, Selected Papers of the 20th International
  Conference on Multiple Criteria Decision Making 2009}}, pages 35--46.
  Springer-Verlag Berlin, Heidelberg, 2011.

\bibitem{Hasenjager2005}
M.~Hasenj{\"a}ger and B.~Sendhoff.
\newblock {Crawling Along the Pareto Front: Tales From the Practice}.
\newblock In {\em The 2005 IEEE Congress on Evolutionary Computation (IEEE CEC
  2005)}, pages 174--181, Piscataway, NJ, 2005. IEEE Press.

\bibitem{HeiMieNie06}
E.~Heikkola, K.~Miettinen, and P.~Nieminen.
\newblock {Multiobjective Optimization of an Ultrasonic Transducer Using
  NIMBUS}.
\newblock {\em Ultrasonics}, 44:368--380, 2006.

\bibitem{Cplex_homepage}
{IBM ILOG CPLEX Optimization Studio Home Page}, Last accessed: \today.
\newblock url =
  \url{http://www-01.ibm.com/software/integration/optimization/cplex-optimizat%
ion-studio/}.

\bibitem{Korhonen1996}
P.~Korhonen and J.~Wallenius.
\newblock {Behavioural Issues in MCDM: Neglected Research Questions}.
\newblock {\em Journal of Multi-Criteria Decision Analysis}, 5:178--182, 1996.

\bibitem{Larichev1992}
O.~I. Larichev.
\newblock {Cognitive Validity in Design of Decision-Aiding Techniques}.
\newblock {\em Journal of Multi-Criteria Decision Analysis}, 3:127 -- 138,
  1992.

\bibitem{Laukkanen2010}
T.~Laukkanen, T.-M. Tveit, V.~Ojalehto, K.~Miettinen, and C.-J. Fogelholm.
\newblock {An Interactive Multi-Objective Approach to Heat Exchanger Network
  Synthesis}.
\newblock {\em Computers \& Chemical Engineering}, 34:943--952, 2010.

\bibitem{Lotov2004}
A.~V. Lotov, V.~A. Bushenkov, and G.~A. Kamenev.
\newblock {\em {Interactive Decision Maps}}.
\newblock Kluwer Academic Publishers, Boston, 2004.

\bibitem{Miettinen1999}
K.~Miettinen.
\newblock {\em {Nonlinear Multiobjective Optimization}}.
\newblock Kluwer Academic Publishers, Boston, 1999.

\bibitem{Miettinen2006b}
K.~Miettinen.
\newblock {IND-NIMBUS for Demanding Interactive Multiobjective Optimization}.
\newblock In T.~Trzaskalik, editor, {\em Multiple Criteria Decision Making'05},
  pages 137--150. The Karol Adamiecki University of Economics in Katowice,
  Katowice, 2006.

\bibitem{Miettinen1995}
K.~Miettinen and M.~M\"akel\"a.
\newblock {Interactive Bundle-based Method for Nondifferentiable Multiobjective
  Optimization: NIMBUS}.
\newblock {\em Optimization}, 34:231--246, 1995.

\bibitem{Miettinen2000}
K.~Miettinen and M.~M\"akel\"a.
\newblock {Interactive Multiobjective Optimization System WWW-NIMBUS on the
  Internet}.
\newblock {\em Computers \& Operations Research}, 27:709--723, 2000.

\bibitem{Miettinen2006}
K.~Miettinen and M.~M. M\"akel\"a.
\newblock {Synchronous Approach in Interactive Multiobjective Optimization}.
\newblock {\em European Journal of Operational Research}, 170:909--922, 2006.

\bibitem{Miettinen2008}
K.~Miettinen, F.~Ruiz, and A.~P. Wierzbicki.
\newblock {Introduction to Multiobjective Optimization: Interactive
  Approaches}.
\newblock In J.~Branke, K.~Deb, K.~Miettinen, and R.~Slowinski, editors, {\em
  {Multiobjective Optimization: Interactive and Evolutionary Approaches}},
  pages 27--57. Springer-Verlag Berlin, Heidelberg, 2008.

\bibitem{Mie07}
Kaisa Miettinen.
\newblock {Using Interactive Multiobjective Optimization in Continuous Casting
  of Steel}.
\newblock {\em Materials and Manufacturing Processes}, 22:585--593, 2007.

\bibitem{Monz2008}
M.~Monz, K.~H. Kufer, T.~R. Bortfeld, and C.~Thieke.
\newblock {Pareto Navigation - Algorithmic Foundation of Interactive
  Multi-Criteria IMRT Planning}.
\newblock {\em Physics in Medicine and Biology}, 53:985--998, 2008.

\bibitem{Nakayama2009}
H.~Nakayama, Y.~Yun, and M.~Yoon.
\newblock {\em {Sequential Approximate Multiobjective Optimization Using
  Computational Intelligence}}.
\newblock Springer-Verlag Berlin, Heidelberg, 2009.

\bibitem{Octave_homepage}
{Octave Home Page}.
\newblock url = \url{www.gnu.org/software/octave/}, Last accessed: \today.

\bibitem{Philips2009}
H.~M. Phillips, K.~E. Sahlstedt, K.~Frank, J.~Bratby, W.~Brennan, S.~Rogowski,
  D.~Pier, W.~Anderson, M.~Mulas, J.~B. Copp, and N.~Shirodkar.
\newblock {Wastewater Treatment Modelling in Practice: a Collaborative
  Discussion of the State of the Art}.
\newblock {\em Water Science and Technology}, 59:695--704, 2009.

\bibitem{Price1983}
W.~L. Price.
\newblock {Global Optimization by Controlled Random Search}.
\newblock {\em Journal of Optimization Theory and Applications}, 40:333--348,
  1983.
\newblock 10.1007/BF00933504.

\bibitem{RuoBomMieTer09}
H.~Ruotsalainen, E.~Boman, K.~Miettinen, and J.~Tervo.
\newblock {Nonlinear Interactive Multiobjective Optimization Method for
  Radiotherapy Treatment Planning with Boltzmann Transport Equation}.
\newblock {\em Contemporary Engineering Sciences}, 2:391--422, 2009.

\bibitem{Ruotsalainen2010}
H.~Ruotsalainen, K.~Miettinen, J.-E. Palmgren, and T.~Lahtinen.
\newblock {Interactive Multiobjective Optimization for Anatomy based
  Three-Dimensional HDR Brachytherapy}.
\newblock {\em Physics in Medicine and Biology}, 55:4703--4719, 2010.

\bibitem{RuzikaWiecek2005}
S.~Ruzika and M.~M. Wiecek.
\newblock {Approximation Methods in Multiobjective Programming}.
\newblock {\em Journal of Optimization Theory and Applications}, 126:473--501,
  2005.

\bibitem{Sahlstedt2010}
K.~Sahlstedt, J.~Hakanen, and K.~Miettinen.
\newblock {Interactive Multiobjective Optimization in Wastewater Treatment
  Plant Operation and Design}.
\newblock In {\em {Proceedings of the ECWATECH-2010, the IWA Specialist
  Conference 'Water and Wastewater Treatment Plants in Towns and Communities of
  the XXI Century: Technologies, Design and Operation', CD-ROM}}, 2010.

\bibitem{Sawaragi1985}
Y.~Sawaragi, H.~Nakayama, and T.~Tanino.
\newblock {\em {Theory of Multiobjective Optimization}}.
\newblock Academic Press, Orlando, 1985.

\bibitem{Henterryck1999}
P.~van Hentenryck.
\newblock {\em {The OPL Optimization Programming Language}}.
\newblock The MIT Press, Cabridge, 1999.

\bibitem{Wierzbicki1980b}
A.~P. Wierzbicki.
\newblock {The Use of Reference Objectives in Multiobjective Optimization}.
\newblock In G.~Fandel and T.~Gal, editors, {\em Multiple Criteria Decision
  Making Theory and Applications}, pages 468--486. Springer-Verlag Berlin,
  Heidelberg, 1980.

\bibitem{Wierzbicki1986}
A.~P. Wierzbicki.
\newblock {On the Completeness and Constructiveness of Parametric
  Characterizations to Vector Optimization Problems}.
\newblock {\em OR Spectrum}, 8:73--87, 1986.

\bibitem{Wilson2001}
B.~Wilson, D.~Cappelleri, T.~W. Simpson, and M.~Frecker.
\newblock {Efficient Pareto Frontier Exploration Using Surrogate
  Approximations}.
\newblock {\em Optimization and Engineering}, 2:31--50, 2001.

\bibitem{Xu2004}
L.~Xu, T.~Reinikainen, W.~Ren, B.~P. Wang, Z.~Han, and D.~Agonafer.
\newblock {A Simulation-Based Multi-Objective Design Optimization of Electronic
  Packages Under Thermal Cycling and Bending}.
\newblock {\em Microelectronics Reliability}, 44:1977 -- 1983, 2004.

\bibitem{Yapicioglu2011}
H.~Yapicioglu, H.~Liu, A.~E. Smith, and G.~Dozier.
\newblock {Hybrid Approach for Pareto Front Expansion in Heuristics}.
\newblock {\em Journal of Operational Research Society}, 62:348--359, 2011.

\end{thebibliography}
\bibliographystyle{plain}
\end{document}